\documentclass[12pt]{article}
\usepackage[utf8]{inputenc}
\usepackage{amsmath}
\usepackage{amsfonts}
\usepackage{amssymb}
\usepackage{amsthm}
\usepackage{relsize}
\usepackage{array}
\usepackage{multirow}
\usepackage{tabu}
\usepackage{soul}
\usepackage{graphicx}
\usepackage{epigraph}
\usepackage{float}
\usepackage[english]{babel}
\usepackage[usenames, dvipsnames]{color}
\usepackage{fancyhdr}
\usepackage[Sonny]{fncychap}
\usepackage{tikz-cd}
\usepackage{tkz-euclide}
\usepackage{sseq}
\usepackage{newfloat}
\usepackage{amsmath,mathtools}
\usepackage{pgfplots}
\usepackage{geometry}
 \usepackage{sectsty}
\sectionfont{\fontsize{14}{15}\selectfont}

\usepackage[all]{xy}
\DeclareFloatingEnvironment[fileext=dia]{diagram}

\theoremstyle{definition}

\newtheorem{remark}{Remark}

\usepackage[most]{tcolorbox}

\makeatletter
\newtcolorbox{note}[1][]{%
	breakable,
	enhanced jigsaw, 
	borderline west={3pt}{0pt}{black!10!white}, 
	borderline south={1pt}{0pt}{black!10!white}, 
	borderline east={1pt}{0pt}{black!10!white},
	borderline north={1pt}{0pt}{black!10!white},
	sharp corners, 
	boxrule=0pt, 
	attach title to upper, 
	left=0pt,
	right=0pt,
	top=0pt,
	bottom=0pt,
	boxsep=5pt,
	colback=white,
	frame hidden,
	#1
}
\makeatother

\makeatletter
\newtcolorbox{note1}[1][]{%
	breakable,
	enhanced jigsaw, 
	sharp corners, 
	boxrule=0pt, 
	attach title to upper, 
	fontupper=\linespread{1.1}\fontfamily{qpl}\selectfont,
	fontlower=\linespread{1.1}\fontfamily{qpl}\selectfont, 
	left=0pt,
	right=0pt,
	top=0pt,
	bottom=0pt,
	boxsep=3pt,
	colback=green!3!white,
	frame hidden,
	before skip=10pt plus 2pt,after skip=10pt plus 2pt,
	#1
}
\makeatother

\makeatletter
\newcommand\tabfill[1]{%
	\dimen@\linewidth
	\advance\dimen@\@totalleftmargin
	\advance\dimen@-\dimen\@curtab
	\parbox[t]\dimen@{#1\ifhmode\strut\fi}%
}
\makeatother

\usepackage{tabularx}
\usepackage{imakeidx}
\usepackage{hyperref}
\hypersetup{colorlinks={true},linkcolor={blue},citecolor={blue},urlcolor={blue}, linktoc=all}  
 
 \usepackage[capitalize, nameinlink]{cleveref}
 \crefdefaultlabelformat{#2\textbf{#1}#3}
 \crefname{figure}{Figure}{Figures} 
 \Crefname{figure}{Figure}{Figures}
 \crefname{table}{Table}{Tables}
 \Crefname{table}{Table}{Tables}
 \crefname{section}{\S\hspace{-1mm}}{\S\hspace{-1mm}}
 \Crefname{section}{\S\hspace{-1mm}}{\S\hspace{-1mm}}
 \crefname{equation}{}{}
 \Crefname{equation}{}{}
 \crefname{example}{Geometric Pattern}{Geometric Patterns} 
 \Crefname{example}{Geometric Pattern}{Geometric Patterns}

 \usepackage{graphicx,tipa}
 \newcommand{\arc}[1]{{%
 		\setbox9=\hbox{#1}%
 		\ooalign{\resizebox{\wd9}{\height}{\texttoptiebar{\phantom{A}}}\cr#1}}}
 
 \usepackage{footnote}

\begin{document}

\title{\textbf{Circular Figures in Elamite Mathematics}}

\author{Nasser Heydari\footnote{Email: nasser.heydari@mun.ca}~ and  Kazuo Muroi\footnote{Email: edubakazuo@ac.auone-net.jp}}

\maketitle

\begin{abstract}
 In this article, we study a particular group of plane figures whose    constants are listed in the Susa Mathematical Tablet No.\,3 (\textbf{SMT No.\,3}).  We explain   possible ways to define these figures and seek to demonstrate that the Susa scribes   used  complicated calculations to obtain  such numbers. We also give examples of circular figures used for    Elamite artifacts  through which one can appreciate the significance of these figures   in   Elamite art.  
\end{abstract}

\section{Introduction}
In Babylonian mathematics  there are mathematical tablets whose texts list  many geometrical coefficients regarding geometric plane figures.  Many scholars call them  \textit{coefficient lists} or \textit{constant lists}.  In these texts, a list of numbers and modifiers are listed together. However, no information regarding the related calculations are provided.  By checking  the possible calculations  required   to find   some of these constants, one might better understand both the complexity of the calculations and  the level of mathematical skill employed by  the  scribes who recorded them.      

One of the famous coefficient lists in Babylonian mathematics is     \textbf{SMT No.\,3}, which   is one of   26 clay tablets excavated from Susa in  southwest Iran by French archaeologists in 1933 and now held in the collection of the Louvre. The texts of all the Susa mathematical tablets (\textbf{SMT}) along with their interpretations were first published in 1961 (see \cite{BR61}).

This  tablet\footnote{The reader can see the new  photos of this tablet on the website of the Louvre's collection. Please see \url{https://collections.louvre.fr/en/ark:/53355/cl010185653} for obverse  and   reverse.}  contains a list of constants with 68 entries, almost all of which are preserved in good condition. The content of this tablet is as follows: line 1 is the headline, lines 2-32 contain the mathematical constants regarding areas and dimensions of geometrical figures  and  lines 33-69 contain the non-mathematical constants  concerning work quotas, metals, and so on. The last two lines (lines 70-71)  are an example of how to use the constants  of work quotas.

Some of the figures mentioned in this text are rather complex geometric figures and finding their area  requires  mathematical skill and experience.  This clearly shows the considerable interest of the Susa scribes in geometry and their  ability to compute the area and other geometric constants of these figures.

\section{Circular Figures}
\begin{figure}[H]
	\centering
	\includegraphics[scale=1]{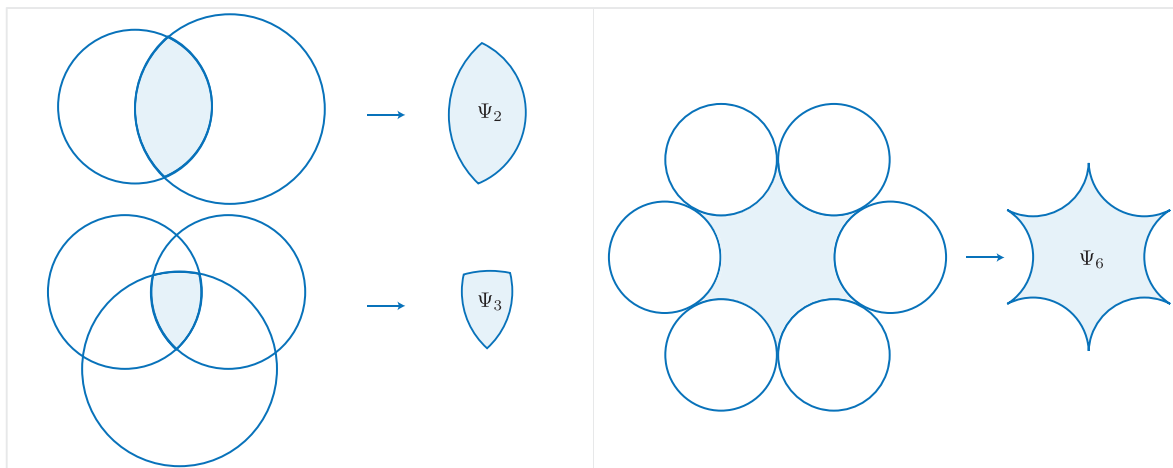}
	\caption{Polyarcs}
	\label{Figure1a}
\end{figure}

Circular figures are usually made up of     arcs of different circles. In other words, these figures can be considered as polygons but their sides instead of line segments are arcs of circles. Unlike polygons, the numbers of sides can be any natural number $n\geq 1$, because we can consider the circle itself as a circular figure with only 1 arc (the whole circumference). For example, intersections     of   circles (not necessarily equal) produce different circular figures. \cref{Figure1a} shows  examples of circular figures with 2, 3 and 6 arcs.   We call the intersection of  two arcs in a circular figure a \textit{vertex}.    If the circular figure has $n$  arcs, we call it  a \textit{polyarc} with $n$ arcs or an  $n$\textit{-arc}    and denote it by   $\Psi_n$.

Although the definition of polyarcs encompasses a wide variety of figures, we are mostly interested in those  where all   arcs are equal and come from the same circle (the \textit{constructive} circle). We prefer to use the term \textit{regular polyarcs} for  such specific  circular figures. For example, the polyarc $\Psi_6$ in \cref{Figure1a} is   regular    while  $\Psi_2$ and $\Psi_3$ are not.  We usually use $\Psi_n(a)$ for a regular polyarc with $n$ arcs  the common length of whose equal arcs   is $a$. Note that this common value is a fraction of the circumference of the constructive circle.  

\begin{figure}[H]
	\centering
	\includegraphics[scale=1]{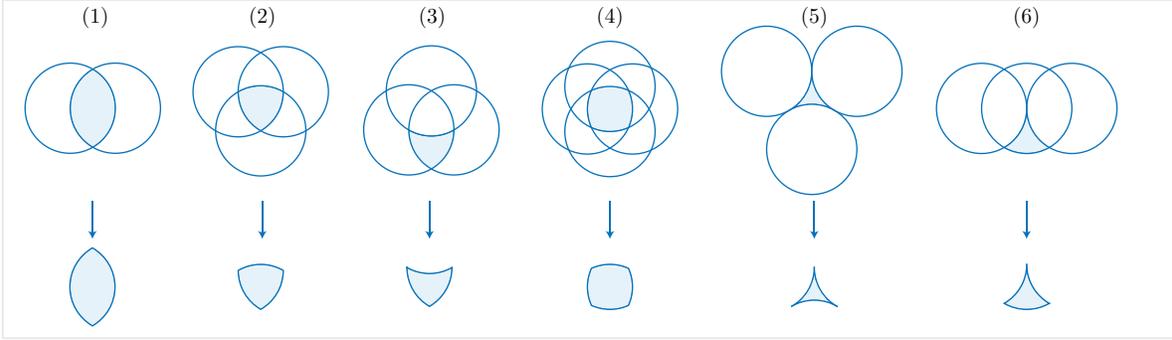}
	\caption{Regular polyarcs}
	\label{Figure1b}
\end{figure}

In this article, we only consider  regular polyarcs whose arcs are all convex or all concave. We call them \textit{regular convex polyarcs} and   \textit{regular concave polyarcs}, respectively.   In \cref{Figure1b}, cases 1, 2 and 4 are regular  convex polyarcs, while case 3 is a regular  concave polyarc. Note that cases 3 and 6 are regular polyarcs but neither is convex nor concave.

An  interesting problem  might be the computation of the area of a regular convex or a  regular concave  polyarc $\Psi_n(a)$. Some of these geometric figures have     properties that   make the accomplishment of this task  possible. In this article, we consider   two groups of regular polyarcs whose areas can be computed because of the way we construct them.  We first discuss the concave case.

\subsection{Regular concave polyarcs}
\begin{figure}[H]
	\centering
	\includegraphics[scale=1]{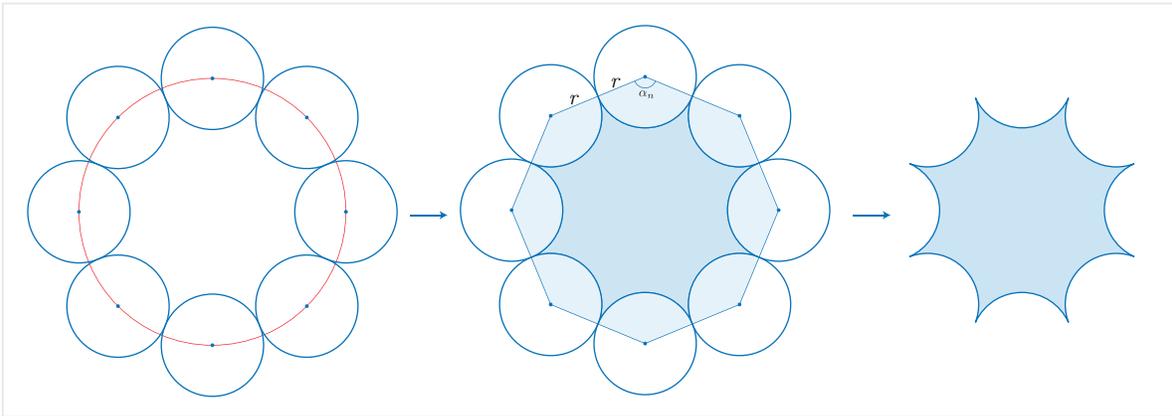}
	\caption{Construction of a   regular concave polyarc}
	\label{Figure1c}
\end{figure}

 One  efficient way to construct this group of regular concave  polyarcs with $n$ arcs   is to consider a   circle  of radius $r$   and repeat (rotate) it $n$ times  along a bigger circle in a way that the obtained $n$ copies   form a chain of pairwise tangent circles.\footnote{It can be shown that the radius of the bigger circle must be $R=\frac{r}{\sin(\frac{\pi}{n})} $.}   If we connect the centers of adjacent circles together, we get a regular $n$-gon whose sides are of length $2r$, where $r$ is the radius of the constructive circle. \cref{Figure1c} shows this situation for the case $n=8$. The symmetry of this figure helps us to obtain a formula for the area of this polyarc with respect to its arc  $a$. We use $\Psi^c_n(a)$ for the  polyarcs in this group.

Let us compute the area of this regular concave  polyarc $\Psi^c_n(a)$. It is clear from \cref{Figure1c} that the area of   $\Psi^c_n(a)$ is the area of the regular $n$-gon with side $2r$ minus $n$ times the area of a sector $\Lambda_n$ of the constructive circle and internal angle $\alpha_n=\frac{(n-2)\pi}{n}$.  Since the circumference is $c=2\pi r$  and  $a$ is the arc length of the internal angle $\alpha_n$,   the  values of circumference $c$ and radius $r$ with respect to $a$  are 
\begin{equation}\label{equ-a}
	c =\frac{2na}{n-2}
\end{equation}   
and
\begin{equation}\label{equ-b}
	r =\frac{na}{(n-2)\pi}.
\end{equation}
It follows from \cref{equ-b} that the area of the sector $\Lambda_n$ is 
\[ S_{\Lambda_n}=\frac{1}{2}\alpha_n r^2= \frac{1}{2}\times \frac{(n-2)\pi}{n} \times \left(\frac{na}{(n-2)\pi}\right)^2 \]
which can be simplified as
\begin{equation}\label{equ-c}
	S_{\Lambda_n}= \frac{na^2}{2\pi(n-2)}.
\end{equation}
On the other hand, the area of the regular $n$-gon  with side $2r$, say $ \Gamma_n(2r) $, is 
\[ S_{\Gamma_n(2r)}= \frac{n}{4}\times   (2r)^2\times \cot\left(\frac{\pi}{n}\right)= nr^2     \cot\left(\frac{\pi}{n}\right). \]
If we use \cref{equ-b} in the last equation, we get
\begin{equation}\label{equ-d}
	S_{\Gamma_n(2r)}= \frac{n^3a^2}{\pi^2(n-2)^2} \times \cot\left(\frac{\pi}{n}\right).
\end{equation}
Now, it follows from \cref{equ-c} and \cref{equ-d} that
\begin{align*}
	S_{\Psi^c_n(a)}& = S_{\Gamma_n(2r)} - n\times S_{\Lambda_n}\\
	                     & =  \frac{n^3a^2}{\pi^2(n-2)^2} \times \cot\left(\frac{\pi}{n}\right) -n\times   \frac{na^2}{2\pi(n-2)}\\
	                     & =  \left(\frac{n^3}{\pi^2(n-2)^2} \times \cot\left(\frac{\pi}{n}\right) -  \frac{n^2}{2\pi(n-2)}\right)a^2.
\end{align*}
Hence, the area of the  regular concave polyarc   $\Psi^c_n(a)$ with $n$ arcs of the same length $a$ is 
\begin{equation}\label{equ-e}
	S_{\Psi^c_n(a)}=  \left(\frac{n^3}{\pi^2(n-2)^2} \times \cot\left(\frac{\pi}{n}\right) -  \frac{n^2}{2\pi(n-2)}\right)a^2.
\end{equation}

\subsection{Regular convex polyarcs}
\begin{figure}[H]
	\centering
	\includegraphics[scale=1]{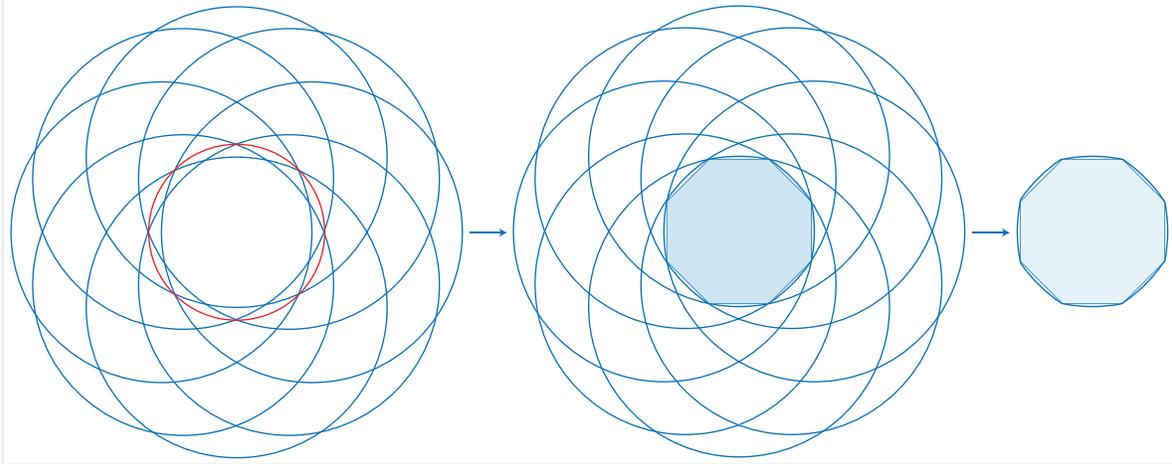}
	\caption{Construction of a regular convex  polyarc}
	\label{Figure1d}
\end{figure}

A  group of regular convex  polyarcs    with $n$ sides denoted by $\Psi^v_n(a)$ can be constructed by repeating (rotating) a constructive circle  $n$ times   around a smaller circle. In this case, the innermost $n$ equal arcs form a regular convex  polyarc. Unlike the concave case, finding the area of such polyarcs is not so easy in general and one may use different approaches such as integration. However, in some special cases ($n=2,3,4,5,6$) one can use the symmetry of the figure and use elementary mathematics to find a formula with respect to the arc length $a$ and $n$.  The key part in this process is to observe that the area of $\Psi^v_n(a)$  is the area  of the inscribed regular $n$-gon plus the total area of $n$ equal circular segments (see \cref{Figure1d}). 

\begin{figure}[H]
	\centering
	\includegraphics[scale=1]{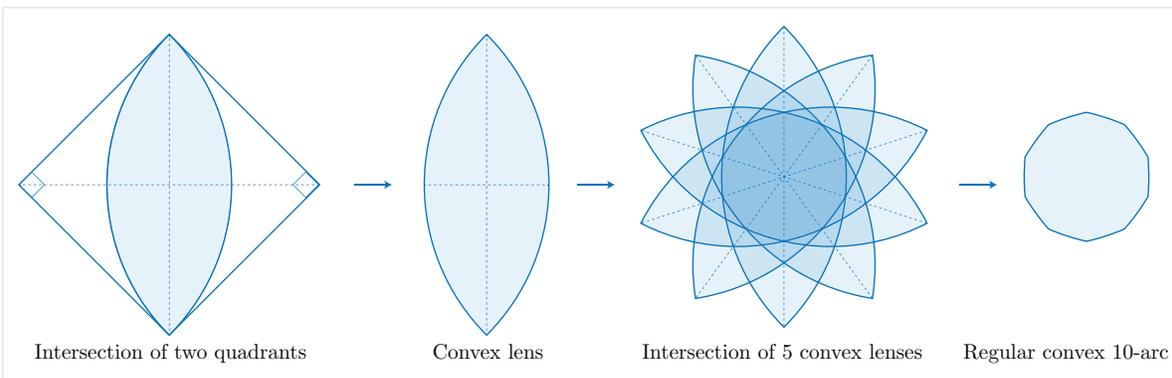}
	\caption{Construction of a regular convex  polyarc with $n=2m$ arcs}
	\label{Figure1e}
\end{figure}

For even numbers $n=2m$, there is an alternate way to define this group of regular convex  polyarcs. In this method, instead of   circles we use a special regular  convex 2-arc. This polyarc is obtained by intersecting   two equal circles such that the length of each arc is the quarter of the circumference of the constructive circle. The obtained polyarc looks like a ``convex lens''.   If we rotate a convex lens $m$ times  around   its center, then the  intersection of the obtained $n$ copies forms a regular convex polyarc with $n=2m$ arcs (see \cref{Figure1e}).  

To compute the area of this regular convex  $n$-arc $\Psi^v_n(a)$, one needs to use the symmetry of the figure and perform some    complicated calculations. We skip the related calculations and only give the formula with respect to $n$ and $r$, the radius of the quadrants:
\begin{equation}\label{equ-f}
	\begin{split}
	S_{\Psi^v_n(a)}&=  n\left[\arcsin\left(\frac{1}{\sqrt{2}}\sin\left(\frac{\pi}{n}\right) \sqrt{1+\cos^2\left(\frac{\pi}{n}\right)}-\frac{1}{2\sqrt{2}}\sin\left(\frac{2\pi}{n}\right)\right)\right]r^2\\
	& - n\left[\frac{1}{2}\sin\left(\frac{\pi}{n}\right)\sqrt{1+\cos^2\left(\frac{\pi}{n}\right)}-\frac{1}{4}\sin\left(\frac{2\pi}{n}\right)\right]r^2.\\
	\end{split}
\end{equation}

\section{Circular figures in \textbf{SMT No.\,3} }
In \textbf{SMT No.\,3}, the constants of several geometric figures are listed some of which have not been fully understood since the publishing of the text of this tablet. The figures we are going to consider here are listed in lines 5,6 and 16-25.  The terms used by the Susa scribe to name these figures might be the reason for this lack of comprehension. There are  six circular figures in these lines which we explain in the following sections.

\subsection{A barley-field}
The transliteration and translation of   lines 16-18 are as follows:
 \begin{note1} 
\underline{Obverse: Lines 16-18} \\
Transliteration:\\
(L16) 13,20 igi-g[ub] \textit{\v{s}\`{a}} a-\v{s}\`{a} \v{s}e\\
(L17) 56,40 dal \textit{\v{s}\`{a}} a-\v{s}\`{a} \v{s}e\\
(L18) 23,20 [\textit{pi}]\textit{-ir-ku \v{s}\`{a}} a-\v{s}\`{a} \v{s}e
 \end{note1}
\noindent
Translation:\\
(L16) 0;13,20 is the constant of a ``barley-field''.\\
(L17) 0;56,40 is the length of a ``barley-field''.\\
(L18) 0;23,20 is the width of a ``barley-field''.\\

\begin{figure}[H]
	\centering
	\includegraphics[scale=1]{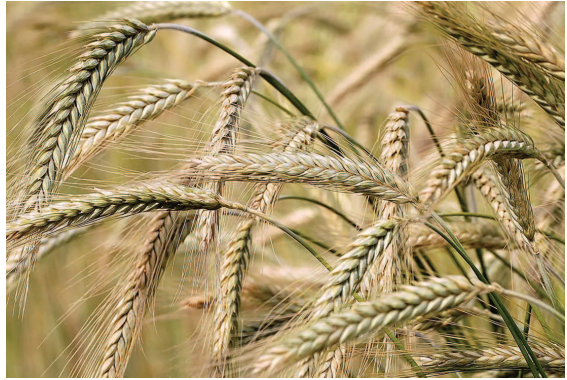}
	\caption{Barley seeds (photo credit: Encyclopædia Britannica)}
	\label{Barley}
\end{figure}

The scribe has used the  term \textit{barley-field} for this geometric figure which might be derived  from the shape of a barley seed (see \cref{Barley}).

\begin{figure}[H]
	\centering
	\includegraphics[scale=1]{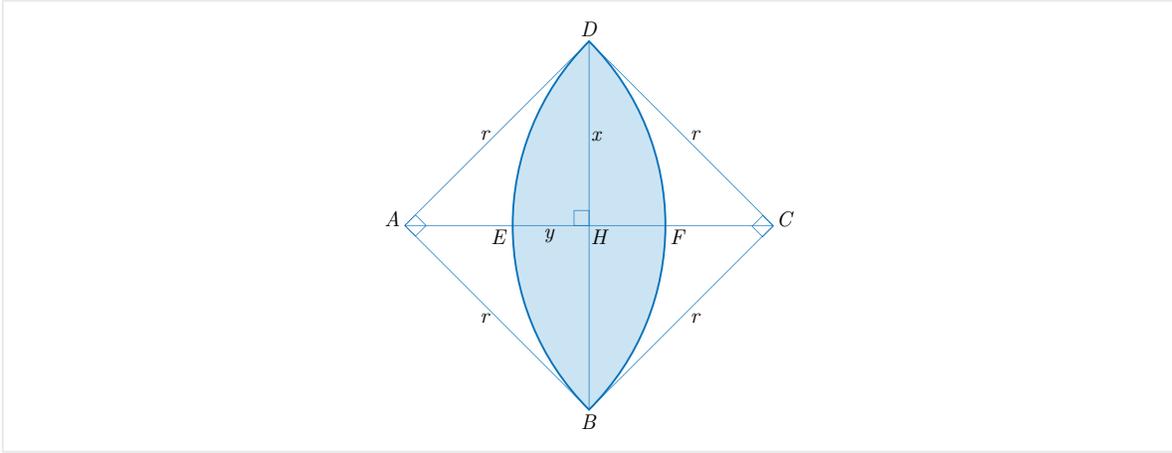}
	\caption{A barley-field figure}
	\label{Figure1}
\end{figure}

 In fact, the figure is   a regular convex 2-arc   (or a convex lens). As mentioned before, it is the intersection of two equal circles such that they separate one-fourth of the circumference of one  another (see \cref{Figure1e}). Beside the figure itself, the Susa scribe has listed two constants for its  length and its width. The former is the line segment connecting  the vertices of the figure and the latter is the line segment connecting the midpoints of the two arcs. We try to verify the scribe's numbers.   To see that, consider a convex lens, say $\Sigma$, whose common arcs have length $a$.  Assume that the radius of the equal quadrants is $r$ and the  length   and the width  are   $2x$ and $2y,$ respectively (see \cref{Figure1}). We   compute the values of $x$, $y$  and  the area of the barley-field figure $\Sigma$ with respect to $a$.

It is clear from  \cref{Figure1} that   $\overline{DB}=2x $ and $\overline{EF}=2y $. Since  the arc  $\arc{DFB} $   is one-fourth of the circumference $c$ of  a circle with radius $r$  and   $c=2\pi r  $ it follows that 
\[ a=\frac{c}{4}=\frac{2\pi r}{4}=\frac{\pi r}{2}  \]
and consequently
\begin{equation}\label{equ-aa}
	r =\frac{2a}{\pi}.	
\end{equation}   

In the right triangle  $\triangle DAB $, two sides $AD$ and $AB$ (radii of the circle) are equal, so two angles $ \angle ADH$ and $ \angle ABH$ are equal and their common measure is $45^{\circ}$. This follows that two angle $ \angle DAH$ and $ \angle  BAH$ have the same measure $45^{\circ}$.    So in  the   right triangle   $\triangle AHD $, we have  $\overline{AH}=\overline{DH}  $. A similar argument shows that $\overline{AH}=\overline{BH}  $, $\overline{BH}=\overline{CH}  $ and $\overline{CH}=\overline{DH}  $. In other words, they all have the same length $x$.  Also note  that by the symmetry we have $\overline{EH}=\overline{FH}=y$.

 Now, by \cref{equ-aa} and   the Pythagorean theorem in the right triangle $\triangle AHD$, we can write

\begin{align*}
	~~&~~\overline{AH}^2+\overline{DH}^2=\overline{AD}^2\\
	\Longrightarrow~~&~~2x^2=r^2 \\
	\Longrightarrow~~&~~x=\frac{\sqrt{2}}{2}r \\
	\Longrightarrow~~&~~x= \frac{\sqrt{2}}{2}  \times\frac{2a}{\pi},
\end{align*}  
and therefore
\begin{equation}\label{equ-ab}
	x=\dfrac{\sqrt{2}a}{\pi}.	
\end{equation}   
On the other hand, since $ \overline{AF}=\overline{AH}+ \overline{FH}$, we obtain $r=x+y $. So, it follows from \cref{equ-ab} that 
\[ y=r-x=\frac{2a}{\pi}-\dfrac{\sqrt{2}a}{\pi}, \]
or equivalently
\begin{equation}\label{equ-ac}
	y=\left(\dfrac{2-\sqrt{2}}{\pi}\right)a.	
\end{equation}

To compute the area of the figure, we should note that it is two times the area of the segment with chord $DB$ and arc $\arc{DFB}$. This area in turn is the area of the  quadrant $ADFB$ with radius $r$ minus that of the triangle $\triangle ABD$. So it follows from \cref{equ-aa}   that
\begin{align*}
	S_{\Sigma}&=2\times \left(S_{ADFB}-S_{\triangle ABD}\right)\\
	&=2\times \left(\frac{1}{4}\pi r^2-\frac{\overline{AB}\times \overline{AD}}{2}\right)\\
	&=2\times \left(\frac{1}{4}\pi r^2-\frac{r^2}{2}\right)\\
	&=  \left(\frac{\pi}{2}  -1\right)r^2\\
	&=  \left(\frac{\pi}{2}  -1\right) \times \left(\frac{2a}{\pi}\right)^2 \\
	&=   \left(\frac{2(\pi-2)}{\pi^2}\right) \times a^2,
\end{align*} 
hence
\begin{equation}\label{equ-ad}
	S_{\Sigma}=  \frac{2(\pi-2)a^2}{\pi^2}.	
\end{equation} 
Note that by setting $n=2$ in \cref{equ-f} and using \cref{equ-aa}, we get the same formula. 

Let us verify the scribe's numbers by assuming $a=1$ and    using  \cref{equ-ab}, \cref{equ-ac} and \cref{equ-ad}. Since in Babylonian mathematics  the approximate values $\pi\approx 3 $ and $\sqrt{2}\approx \frac{17}{12}=1;25 $  were usually used in  calculations, we have
\[\overline{DB}=2x=  \dfrac{2\sqrt{2}a}{\pi}\approx \dfrac{2\times (1;25)\times 1}{3} = (2;50)\times (0;20)=0;56,40  \] 
which is the number listed in line 17. On the other hand, for the width of $\Sigma$, we have
\begin{align*}
	\overline{EF} =2y =  2\left(\dfrac{2-\sqrt{2}}{\pi}\right)a \approx  \frac{2}{3}\times(2-1;25 )\times 1=(0;40)\times (0;35)   = 0;23,20
\end{align*} 
which is the number in line 18. Finally, for approximate  area  of the  barley-field  $\Sigma: DEBF$, we can use \cref{equ-ad} to write
\begin{align*}
	S_{\Sigma} = \frac{2(\pi-2)a^2}{\pi^2}  \approx \frac{2(3-2)\times 1^2}{3^2} = \frac{2}{9} = 0;13,20
\end{align*}
which is the   number   in line 16. So,  it seems that  the Susa scribe was capable of performing such calculations  to obtain these numbers.

\subsection{An ox-eye}
The transliteration and translation of   lines 19-21 are as follows:
\begin{note1} 
\underline{Obverse: Lines 16-18} \\
Transliteration:\\
(L19) 16,52,30 igi-gub \textit{\v{s}\`{a}} igi-gud\\
(L20) 52,30 dal \textit{\v{s}\`{a}} igi-gud\\
(L21) 30 \textit{pi-ir-ku \v{s}\`{a}} igi-gud
\end{note1}
\noindent
Translation:\\
(L19) 0;16,52,30 is the constant of an ``ox-eye''.\\
(L20) 0;52,30 is the length of an ``ox-eye''.\\
(L21) 0;30 is the width  of an ``ox-eye''.\\

\begin{figure}[H]
	\centering
	\includegraphics[scale=1]{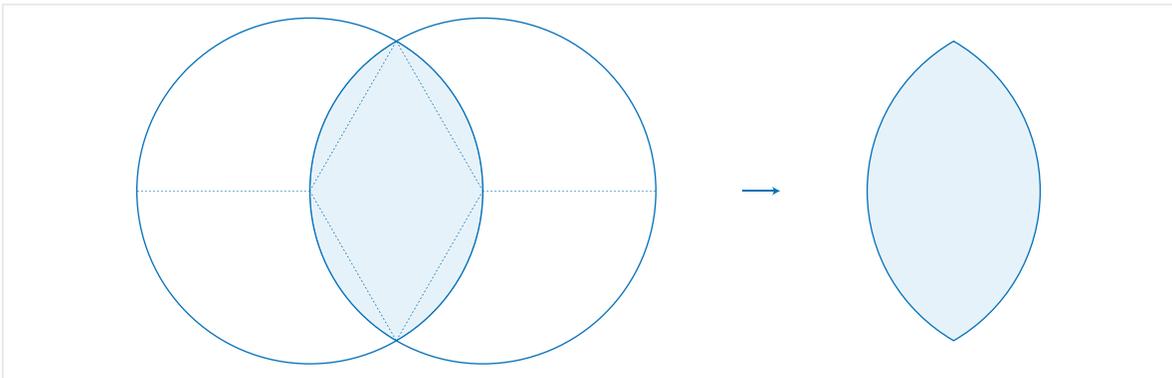}
	\caption{An ox-eye}
	\label{Figure2}
\end{figure}

The    term \textit{ox-eye} for this geometric figure clearly refers to the shape of the eye of an ox. The reason for this naming might be that the figure    is similar to  a convex lens but  it is wider.  It can be made by intersecting   two equal circles such that their arcs are one-third of the circumference of the circle (see \cref{Figure2}).

\begin{figure}[H]
	\centering
	\includegraphics[scale=1]{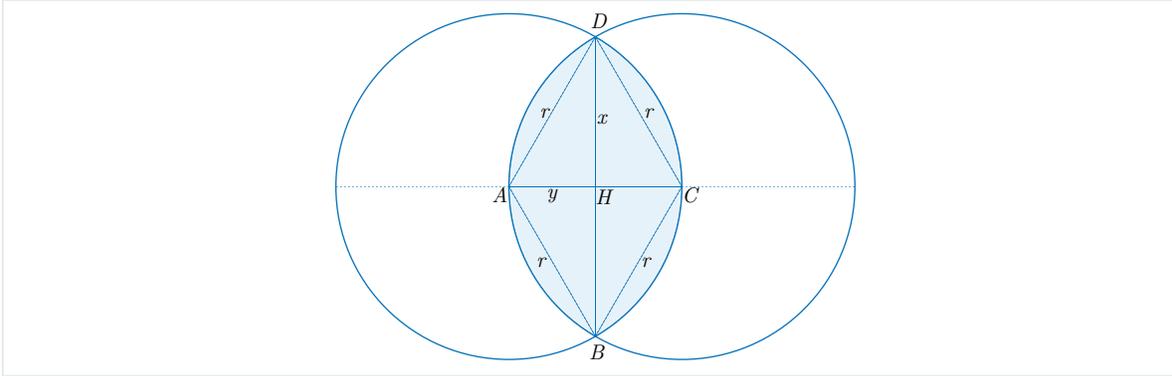}
	\caption{Area of an ox-eye}
	\label{Figure3}
\end{figure}

Three numbers for the area, the  length and the width of the figure has been listed in lines 19-21.     To verify these numbers, consider   the regular convex 2-arc $\Pi$ with arc length $a$ as shown in \cref{Figure3}.   Set $\overline{DH}=x$ and $\overline{AH}=y$.  All sides $AD $, $AB$, $BC$, $DC$ and $AC$ are radii of equal circles, so   the triangles $\triangle ADC $ and $\triangle ABC$ are equilateral and the measures of their internal angles are $60^{\circ} $. Also, the central angles of two arcs   $\arc{DAB} $ and $\arc{DCB} $ are $120^{\circ} $ and the common length $a$ of two arcs   $\arc{DAB} $ and $\arc{DCB} $ is the one-third of the circumference of the circle    with radius  $r$. So $c=2\pi r=3a$ and 
\begin{equation}\label{equ-3-ga}
	r =\frac{3a}{2\pi}.	
\end{equation}
Since $DB$   perpendicularly bisects $\angle ADC $,  the height $DH$ bisects $AC$ and thus $r=\overline{AC}=2\overline{AH}=2y  $. So it follows from \cref{equ-3-ga} that
\begin{equation}\label{equ-3-gb}
	y =\frac{3a}{4\pi}.	
\end{equation}

On the other hand, by the Pythagorean theorem in the right triangle $\triangle AHD $ and using \cref{equ-3-ga} and \cref{equ-3-gb}, we have
\begin{align*}
	x&=\sqrt{r^2-y^2}\\
	&=\sqrt{r^2-\frac{r^2}{4}}\\
	&=\sqrt{\frac{3r^2}{4}}\\
	&=\dfrac{\sqrt{3}}{2}r\\
	&=\dfrac{\sqrt{3}}{2}\times \frac{3a}{2\pi},
\end{align*}
which  implies that
\begin{equation}\label{equ-3-gc}
	x=\frac{3\sqrt{3}a}{4\pi}.	
\end{equation}
Since $\overline{AC}=2y$ and $\overline{BD}=2x$, it follows from \cref{equ-3-ga}, \cref{equ-3-gc}   that
\begin{equation}\label{equ-3-gca}
	\begin{dcases}
		\overline{AC}=\frac{3a}{2\pi},\\
		\overline{BD}=\frac{3\sqrt{3}a}{2\pi}.
	\end{dcases}
\end{equation}

The area  of the ox-eye  $\Pi: ABCD $   is double    the difference between the area  of circular sector  $\Gamma_1: ABCD$ and  the area  of the isosceles triangle $\triangle ABD $. Since the area of the sector $\Gamma_1: ABCD$ is the one-third of the area of a circle with radius $r$, it follows from \cref{equ-3-ga}, \cref{equ-3-gb} and \cref{equ-3-gc}   that
\begin{align*}
	S_{\Pi}&=2\times \left(S_{\Gamma_1}-S_{\triangle ABD}\right)\\
	& =2\times \left[\left(\frac{1}{3}\times \pi r^2\right)-\frac{2xy}{2}\right]\\
	&\approx 2\times \left[ \left(\frac{1}{3}  \times \frac{9\pi a^2}{4\pi^2}\right)   - \left(\frac{3\sqrt{3}a}{4\pi}\times \frac{3a}{4\pi}\right) \right]\\
	&=2\times \left[\frac{3\pi a^2}{4\pi^2}   -  \frac{9\sqrt{3}a^2}{16\pi^2} \right]\\
	&=2 \times\left(\frac{(12\pi-9\sqrt{3}) a^2}{16\pi^2} \right)\\
	& = \frac{(12\pi-9\sqrt{3}) a^2}{8\pi^2},
\end{align*}
or equivalently
\begin{equation}\label{equ-3-gcb}
	S_{\Pi}= \left(\dfrac{3}{2\pi}-\frac{9\sqrt{3}}{8\pi^2}\right) a^2.
\end{equation}

Now, set $a=1$ and use the Babylonian approximations $ \sqrt{3} \approx \frac{7}{4}$ and $ \pi \approx 3 $. It follows from \cref{equ-3-gca} that
\[ \overline{AC}=\frac{3a}{2\pi}\approx \frac{3\times 1}{2\times 3}=\frac{1}{2}=0;30, \]
which is listed in line 21. Also
\[ \overline{BD}=\frac{3\sqrt{3}a}{2\pi}\approx \frac{3\times (7/4)\times 1}{2\times 3}=\frac{7}{8}=0;52,30. \]
which is written in line 20.  Finally, the approximate area is obtained from \cref{equ-3-gcb} as follows:
\[ S_{\Pi}= \left(\dfrac{3}{2\times 3}-\frac{9\times(7/4)}{8\times 9}\right) \times 1\approx   \frac{1}{2} - \frac{7}{32}=\frac{9}{32}=0;16,52,30. \]
This is the very number in line 19 which is given as the constant   of the ``ox-eye''. Thus, the numbers of the Susa scribe are verified.

\subsection{Overlap of two barley-fields}
The transliteration and translation of   line  5 are as follows:
\begin{note1} 
	\underline{Obverse: Line 5} \\
	Transliteration:\\
(L5)   16(sic) [igi]-gub \textit{\v{s}\`{a}} g\'{u}r \textit{\v{s}\`{a}} 2 \v{s}e \textit{i-na} \v{s}\`{a} g\'{u}r gar
\end{note1}
\noindent
Translation:\\
(L5)  0;16  is the constant of a circular figure  with two barley figures  in the middle.\\

\begin{figure}[H]
	\centering
	\includegraphics[scale=1]{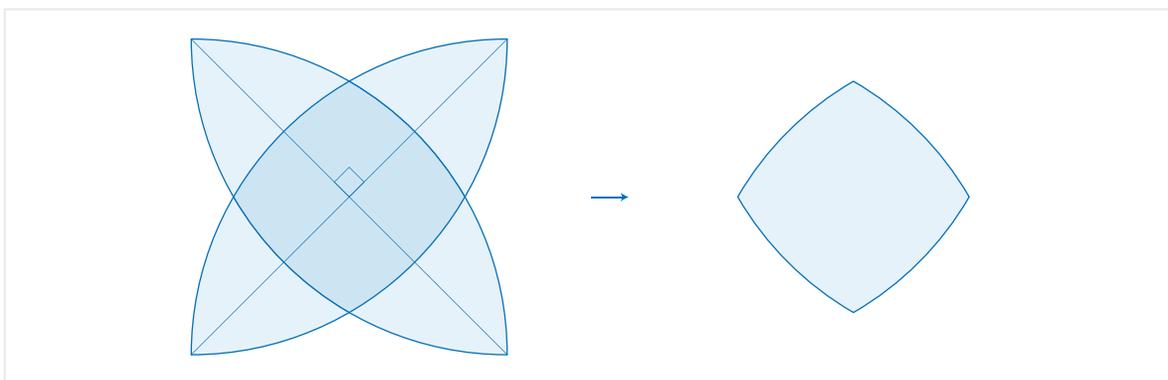}
	\caption{Overlap of two barley-fields}
	\label{Figure4}
\end{figure}

 In  this line, a circular figure  is described as the overlap of  two  barley-field figures. This is exactly the intersection of two convex lenses such that their centers   coincide and their lengths (diagonals) form a right angle. So, it is  a regular convex polyarc with 4 arcs (see \cref{Figure4}).

\begin{figure}[H]
	\centering
	\includegraphics[scale=1]{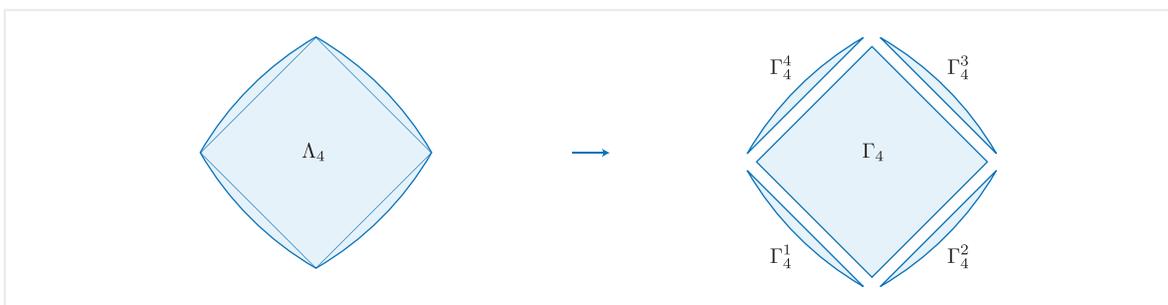}
	\caption{Partition of a regular convex  4-arc}
	\label{Figure5}
\end{figure}

The number 0;16 is written as the constant of this regular polyarc in line 5.  Apparently, the Susa scribe has assumed the common radius of all quadrants forming the figure to be 1.   Let us compute the area of the figure  for a general case when the radius of all quadrants is   $r$. We denote this polyarc  by $\Lambda_4 $. As   in      \cref{Figure5}, this  figure is the union of a quadrilateral $\Gamma_4$ and four   circular segments $\Gamma_4^1, \Gamma_4^2,\Gamma_4^3$  and $\Gamma_4^4 $.  

\begin{figure}[H]
	\centering
	\includegraphics[scale=1]{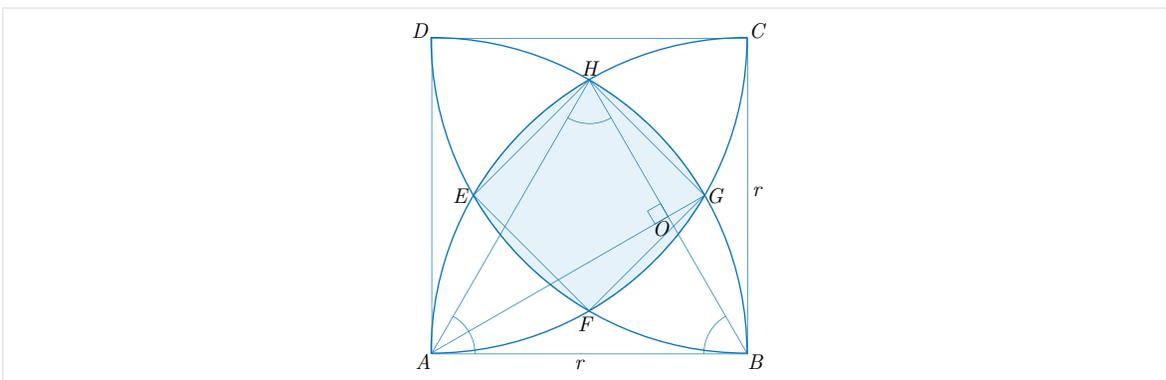}
	\caption{Partition of a  regular convex 4-arc}
	\label{Figure6}
\end{figure}

Now, consider \cref{Figure6} and connect vertices  $H,G$ to the vertex  $A$ to get the isosceles triangle  $\triangle AGH$. Connect $H$ to $B$. Since three line segments $AH$, $AB$ and $HB$ are the radius of equal quadrants, the triangle $\triangle AHB $ is equilateral and all of  its internal angles are $60$ degrees. This implies that  
\[  \angle DAH=\angle  DAB-\angle  HAB =90^{\circ}-60^{\circ}=30^{\circ}. \]  
Because of symmetry, the same is true for angles $\angle GAB$ and $\angle HAG $. Therefore, three arcs $\arc{DH}$, $\arc{HG}$ and $\arc{GB}$ have the same length which is one-twelfth of the circumference of a circle with radius $r$.  Similar arguments can be used to show that all other arcs   $\arc{GF}$, $\arc{FE}$ and $\arc{EH}$ have the common length of $\arc{HG}$, so they all    are equal implying that  their corresponding chords are equal too. Thus,  the sides of the quadrilateral $\Gamma_4:  EFGH $ have the same length.

Next, consider  the isosceles triangle $\triangle AGH$ in which $\angle HAG=30^{\circ} $. So  the other two inner angles are  
$$\angle AGH= \angle AHG=\frac{180^{\circ}-30^{\circ}}{2} =75^{\circ} $$ 
which implies that 
$$\angle BHG= \angle AHG - \angle AHO=75^{\circ}-60^{\circ}=15^{\circ}. $$ 
Similarly, we get $\angle EHA=  15^{\circ} $ and thus 
$$ \angle EHG=\angle EHA + \angle AHB+\angle BHG =15^{\circ}+60^{\circ}+15^{\circ}=90^{\circ}.$$ 
By symmetry, we can see that all other three inner angles of the quadrilateral $\Gamma_4(r):  EFGH $ are also right-angled. So far, we have shown that   the quadrilateral $\Gamma_4(r):  EFGH $ has equal sides and equal right angles, therefore it   must be a square.

On the other hand, since $G$ is in the middle of the arc $ \arc{HGB}$, the line segment $AG$ perpendicularly bisects   its chord $HB$. So   $\angle HOG=90^{\circ} $  and since $ \overline{HB}=r$,    we get $ \overline{HO}=\frac{\overline{HB}}{2}= \frac{r}{2}$. By the Pythagorean theorem in the right triangle $\triangle HOA $, we have  
\[ \overline{AO}=\sqrt{\overline{HA}^2-\overline{HO}^2}=\sqrt{r^2-\frac{r^2}{4}}=\frac{\sqrt{3}}{2}r. \]
Since $ \overline{GO}=\overline{AG}-\overline{AO}$, we get
\begin{equation}\label{equ-3-ba}
	\begin{dcases}
		\overline{HO}= \dfrac{r}{2}\\
		\overline{GO}=\left(1-\dfrac{\sqrt{3}}{2}\right)r.
	\end{dcases}
\end{equation}
Note that  the area of the  square  $\Gamma_4 :  EFGH $ in   \cref{Figure6}  is equal to $S_{\Gamma_4}=\overline{HG}^2  $, so we can apply \cref{equ-3-ba} and   the Pythagorean  theorem  in the right triangle $\triangle HOG $  to obtain $\overline{HG} $  as follows:
\begin{align*}
	\overline{HG}  = \sqrt{\overline{HO}^2+\overline{GO}^2 } = \sqrt{\left(\frac{r}{2}\right)^2+\left(r-\frac{\sqrt{3}r}{2}\right)^2}   = \left(\sqrt{2-\sqrt{3}}\right)  r.
\end{align*} 
Thus,
\begin{equation}\label{equ-3-bc}
	\overline{HG} =\left(\sqrt{2-\sqrt{3}}\right)  r, 
\end{equation}
and   the  area of the inner square $\Gamma_4:  EFGH $  is
\begin{equation}\label{equ-3-b}
	S_{\Gamma_4}=(2-\sqrt{3})r^2. 
\end{equation}

Let us    compute the exact area of the convex polyarc $\Lambda_4$. To do so, we first need to compute the areas of the   circular segments  $\Gamma_4^1, \Gamma_4^2,\Gamma_4^3$  and $\Gamma_4^4 $.   It suffices to compute the area of one of them, say $\Gamma_4^3$, because as we saw above they have the same arcs and chords. It is clear from \cref{Figure6} that $ S_{\Gamma_4^3}  $ is the area of the one-twelfth of the circle with radius $r$ minus the area of the isosceles triangle $\triangle AHG $. So, we can write
\begin{align*}
	S_{\Gamma_4^3}&=\frac{\pi r^2}{12}- S_{\triangle AHG}\\
	&=\frac{\pi r^2}{12}- \dfrac{1}{2}\times \overline{HO}\times \overline{AG}\\ 
	&=\frac{\pi r^2}{12}- \dfrac{1}{2}\times \dfrac{r}{2}\times r\\ 
	&=\frac{\pi r^2}{12}- \dfrac{r^2}{4},
\end{align*}
so
\begin{equation}\label{equ-3-bd}
	S_{\Gamma_4^3}=\left(\frac{\pi}{3}- 1\right) \dfrac{r^2}{4}. 
\end{equation}

Now, we can use  \cref{equ-3-b} and \cref{equ-3-bd} to compute the exact area of the convex polyarc $\Lambda_4$ as follows:
\begin{align*}
	S_{\Lambda_4}&=S_{\Gamma_4}+4\times S_{\Gamma_4^3}\\
	&=  (2-\sqrt{3})r^2+4\times  \left(\frac{\pi}{3}- 1\right) \times \dfrac{r^2}{4}\\
	&=(2-\sqrt{3})r^2+\left(\frac{\pi}{3}- 1\right)  r^2\\
	&= \left(\frac{\pi}{3}+1-\sqrt{3}\right)  r^2.	
\end{align*}
Thus, we obtain
\begin{equation}\label{equ-3-be}
	S_{\Lambda_4}= \left(\frac{\pi}{3}+1-\sqrt{3}\right)  r^2.	
\end{equation}
Note that if we set $n=2$ in the general formula \cref{equ-f}, we get \cref{equ-3-be}.

Finally, let us verify the scribe's number. According to Babylonian tradition,  if we set $r=1$ and use the approximation $ \sqrt{3}\approx \frac{7}{4}$ in \cref{equ-3-b}, the approximate area  of $\Gamma_4 $  is
\[ S_{\Gamma_4}=(2-\sqrt{3})\times 1^2\approx 2-\frac{7}{4}=2-1;45=0;15. \]
By comparing this value with the value $0;16$ given in line 5, one might get the impression  that  the Susa scribe     has made a  mistake here, but there may be something hidden here. If we   only set $r=1$   in \cref{equ-3-b} and assume that  $S_{\Gamma_4}\approx 0;16$ then we get 
\begin{align*}
	~~&~~2-\sqrt{3}  \approx 0;16\\
	\Longrightarrow~~&~~ \sqrt{3}  \approx 2-\dfrac{16}{60}\\
	\Longrightarrow~~&~~\sqrt{3}  \approx 2-\dfrac{4}{15}\\
	\Longrightarrow~~&~~\sqrt{3}  \approx  \dfrac{30-4}{15}\\
	\Longrightarrow~~&~~\sqrt{3}  \approx  \dfrac{26}{15}.
\end{align*}  

In other words, the scribe of this tablet   has considered the less common approximation $\sqrt{3}  \approx  \frac{26}{15} $ to approximate the area of the convex polyarc $ \Lambda_4$. Although this approximate value seems to be rare, there are other mathematical texts using  this approximation in their calculations. One specific example is given in the Achaemenid mathematical tablet   \textbf{W 23291}\footnote{This mathematical tablet  belongs to the Late Babylonian period (ca. the first millennium BC). It contains around 26 problems regarding the areas and dimensions of geometric figures like equilateral triangles.} (see \cite{Fri07-2}, for details). It should be noted that this more accurate approximation to $ \sqrt{3} $ can be obtained by using the usual Babylonian method or the approximate formula  $\sqrt{a^2\pm b}\approx a\pm \frac{b}{2a}$ by setting $a=\frac{5}{3}$ and $b=\frac{2}{9}$: 
\begin{align*}
	\sqrt{3}&=\sqrt{\frac{25}{9}+\dfrac{2}{9}}\\
	&=\sqrt{\left(\frac{5}{3}\right)^2+\dfrac{2}{9}}\\
	&\approx \frac{5}{3}+\frac{\frac{2}{9}}{2\times (\frac{5}{3})}\\
	&=\frac{5}{3}+\frac{1}{15}\\
	&=\frac{26}{15}.	
\end{align*}
Therefore, the approximate area of the inscribed square $\Gamma_4$ is    
\begin{equation}\label{equ-3-ab}
	S_{\Gamma_4}=0;16 
\end{equation}
which is the number in line 5.

It is interesting that if we   use the Babylonian approximations  $\pi\approx 3  $  and   $\sqrt{3}\approx\frac{26}{15} $  in \cref{equ-3-be},  the approximate area of the convex polyarc $\Lambda_4$ is
\[S_{\Lambda_4}\approx\left(\dfrac{3}{3}+1-\dfrac{26}{15}\right) \times 1^2 =   \dfrac{4}{15}=0;16\]
which is the approximate area of the inscribed square  $\Gamma_4$ given in \cref{equ-3-ab}. Therefore, the scribe of this tablet might have hit two birds with one stone!

\subsection{Overlap of three barley-fields}
The transliteration and translation of   line  6 are as follows:
\begin{note1} 
	\underline{Obverse: Line 6} \\
	Transliteration:\\
	(L6)   16,26,46,40 igi-gub \textit{\v{s}\`{a}} g\'{u}r \textit{\v{s}\`{a}} 3 \v{s}e \textit{i-na} \v{s}\`{a} g\'{u}r gar 
\end{note1}
\noindent
Translation:\\
(L6)  0;16,26,46,40 is the constant of a circular figure  with three barley figures in the middle.\\

\begin{figure}[H]
	\centering
	\includegraphics[scale=1]{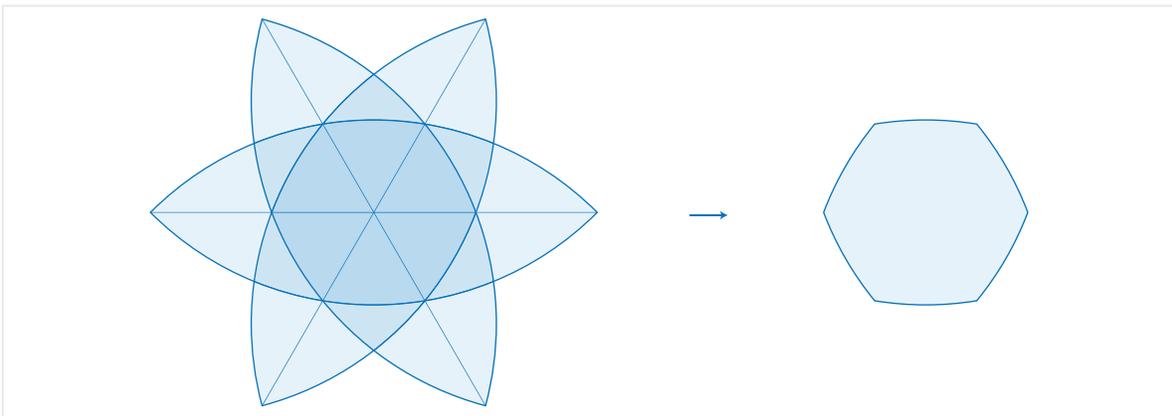}
	\caption{Overlap of three barley-fields}
	\label{Figure7}
\end{figure}

The  scribe refers to the circular figure in this line   as the overlap of  three  barley-field figures. In fact, it is the intersection of three  rotated convex lenses which produces a regular convex polyarc with 6 arcs (see \cref{Figure7}).  The arrangement of barley figures  are  symmetrical meaning that each one is obtained from another   if we rotate them by $60^{\circ}$.

\begin{figure}[H]
	\centering
	\includegraphics[scale=1]{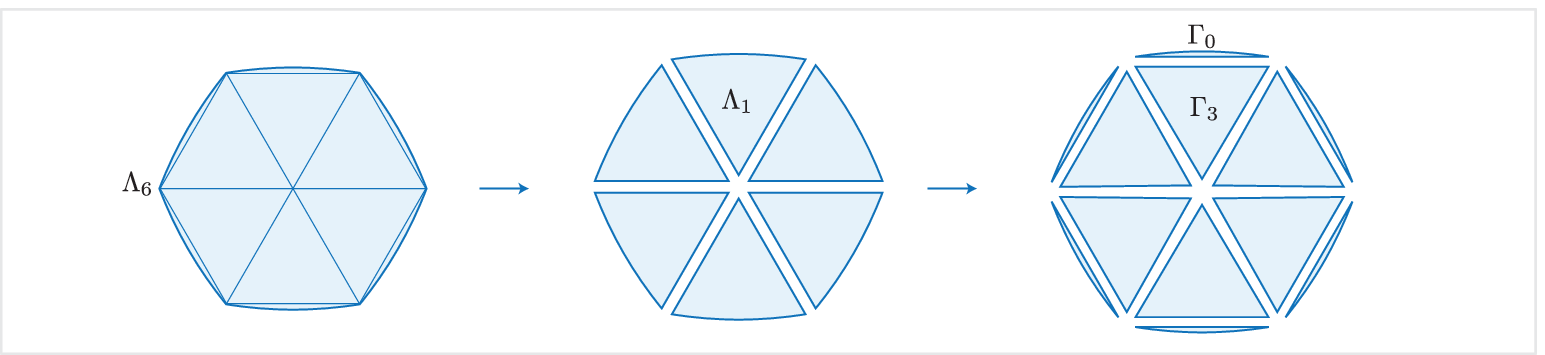}
	\caption{Partition of a  regular convex 6-arc}
	\label{Figure8}
\end{figure}

In line 6, the area of this circular figure, say $ \Lambda_6 $,  is assumed to be $S_{\Lambda_6}=  0;16,26,46,40$.    To calculate the area,   the Susa scribe  of this tablet seems to have assumed the radius  of each quadrant  to be 1   and has omitted the areas of the six small circular segments from his calculations. To verify this, we first compute the area of this regular  convex   polyarc for the general case in which the radius of each quadrant is $r$. If we connect each vertex of the regular polyarc  first  to adjacent vertices and then to its center, we get a partition of the circular figure  consisting  of six   equilateral triangles  and six   small circular segments as shown in  \cref{Figure8}.

In line 6, the area of this circular figure $ \Lambda_6 $  is assumed to be $S_{\Lambda_6}=  0;16,26,46,40$.    To calculate the area,   the Susa scribe  of this tablet seems to have assumed the radius  of each quadrant  to be 1   and has omitted the areas of  six small circular segments from his calculations. To verify this, we first compute the area of this convex hexagon for the general case in which the radius of each quadrant is $r$.

\begin{figure}[H]
	\centering
	\includegraphics[scale=1]{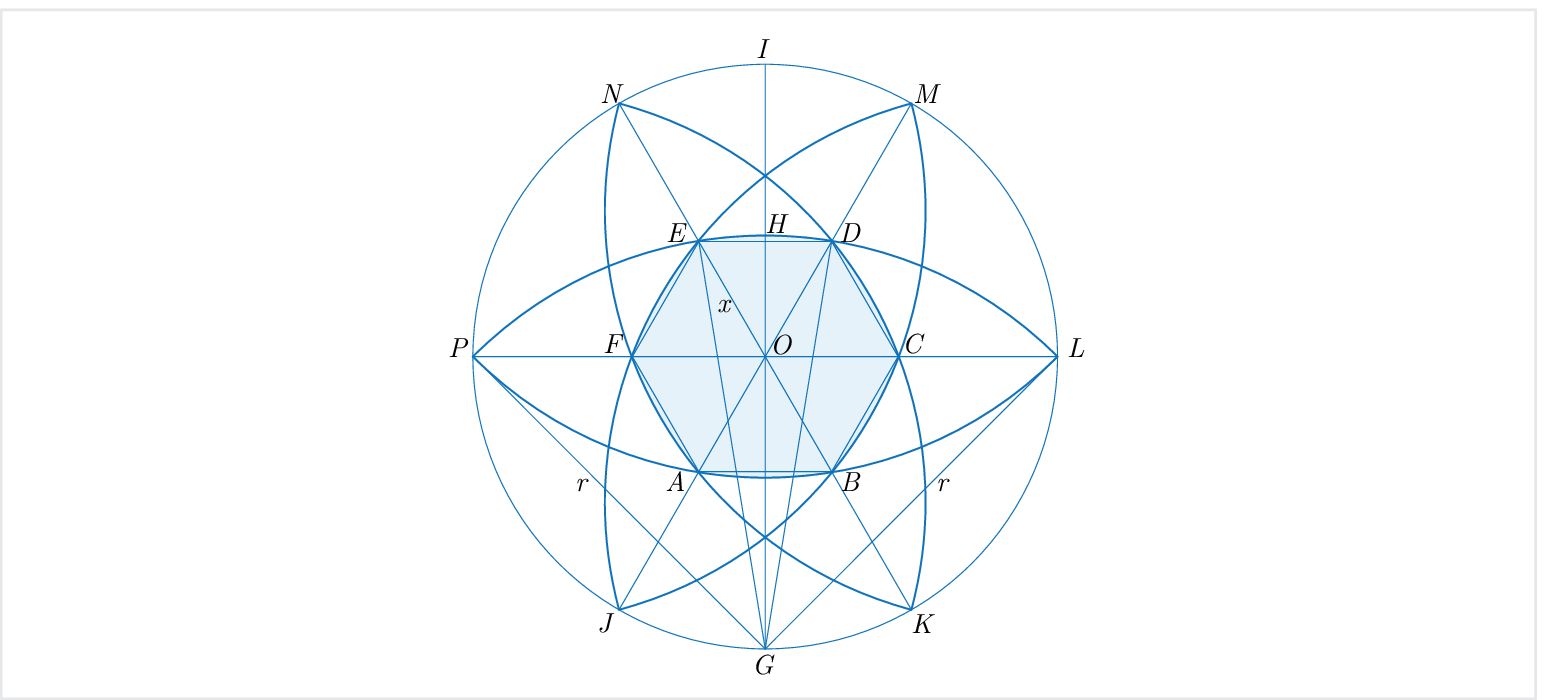}
	\caption{Partition of a regular convex  6-arc}
	\label{Figure9}
\end{figure}

According to  \cref{Figure8}, this convex polyarc is the union of six triangles and six circular segments. We first show that all these triangles are equal and equilateral and also that all the  segments are   equal. In such a case, to find the area we just need to compute the area of the equilateral triangle $\Gamma_3$ and the segment $\Gamma_0$. To do so, we need to consider some auxiliary lines as shown in \cref{Figure9}.

By our assumption, each barley figure  (convex lens) is obtained from the other one by a 60-degree rotation. So, all central angles of the hexagon $\Gamma_6=ABCDEF$ are equal and have the same measure $60^{\circ}$. By symmetry, all line segments connecting $O$ to the vertices of $\Gamma_6$ are equal. Therefore, the hexagon $\Gamma_6$ is regular and all interior triangles partitioning $\Gamma_6$ are equilateral. Let the common length of all sides of all these triangles be $x$. Since $GI$ and $PL$ are perpendicular, $\angle HOE=90^{\circ}-60^{\circ}=30^{\circ}$ which implies that  $\angle OHE=90^{\circ} $. It follows   that in the equilateral triangle $ \triangle OED$, the line segment $OH$ bisects $\angle EOD$. Thus, $OH$ also bisects the base $ED$ and so
\begin{equation}\label{equ-3-ea}
	\overline{EH} = \frac{\overline{ED}}{2}=\frac{x}{2}.	
\end{equation}
Now, by the Pythagorean theorem in the right triangle $\triangle OHE $, we have
\[ \overline{HO} = \sqrt{\overline{OE}^2- \overline{EH}^2}=\sqrt{x^2-\frac{x^2}{4}}=\frac{\sqrt{3}}{2}x.\] 
Since $ \overline{GP}= \overline{GL}=r$, the triangle $ \triangle PGL$  is isosceles and thus the height $GO$ bisects the base $PL$  and we have  $\overline{OK}=\overline{OL} $. By the Pythagorean theorem in the right triangle  $\triangle POG $ we get
\[ \overline{OG}=\sqrt{\overline{OG}^2}=\sqrt{\dfrac{\overline{OP}^2+\overline{OG}^2}{2}} =\sqrt{\dfrac{\overline{PG}^2}{2}} =\dfrac{r}{\sqrt{2}}= \frac{\sqrt{2}}{2}r.\]
It follows from these calculations that
\begin{equation}\label{equ-3-eb}
	\overline{GH}=\frac{\sqrt{3}}{2}x+\frac{\sqrt{2}}{2}r.       	
\end{equation}

Next, by the Pythagorean theorem in the right triangle $ \triangle EHG $ and equations \cref{equ-3-ea} and \cref{equ-3-eb}, we obtain the following quadratic equation  
\[\left(\frac{x}{2}\right)^2+\left(\frac{\sqrt{3}x}{2}+\frac{\sqrt{2}}{2}r\right)^2 = r^2,\]
which can be simplified and written as the standard form  
\begin{equation}\label{equ-3-ec}
	x^2+\left(\frac{\sqrt{6}}{2}r\right)x = \frac{r^2}{2}.
\end{equation} 
Note that $x$ is unknown and $r$ is a constant.

By completing the square, which is a standard method in Babylonian mathematics and called \textit{tak\={\i}ltum} in Akkadian, the solution of the equation  \cref{equ-3-ec}  is obtained as follows:
\begin{align*}
	~~&~~x^2+\left(\frac{\sqrt{6}}{2}\right)x  = \frac{r^2}{2}\\
	\Longrightarrow~~&~~x^2+2\times \left(\frac{\sqrt{6}}{4}r\right)x= \frac{r^2}{2}\\
	\Longrightarrow~~&~~x^2+2\times \left(\frac{\sqrt{6}}{4}r\right)x+\left(\frac{\sqrt{6}}{4}r\right)^2  = \frac{r^2}{2}+\left(\frac{\sqrt{6}}{4}r\right)^2\\
	\Longrightarrow~~&~~\left(x+\frac{\sqrt{6}}{4}\right)^2  = \frac{14}{16}r^2\\
	\Longrightarrow~~&~~x+\frac{\sqrt{6}}{4}r  = \frac{\sqrt{14}}{4}r, 
\end{align*}  
which implies that
\begin{equation}\label{equ-3-f}
	x =\left(\frac{\sqrt{14}-\sqrt{6}}{4}\right)r. 
\end{equation}
It follows from    equation  \cref{equ-3-f}    that the area  of the equilateral triangle   $\triangle OED$  with side $x$ is 
\begin{align*}
	S_{\triangle OED}&= \frac{\sqrt{3}}{4}x^2\\
	&=\frac{\sqrt{3}}{4} \times \left(\frac{\sqrt{14}-\sqrt{6}}{4}\right)^2r^2\\
	&=\frac{\sqrt{3}}{4} \times \left(\frac{14+6-2\sqrt{4\times21}}{16}\right)r^2\\
	&=\frac{\sqrt{3}}{4} \times \left(\frac{5-\sqrt{21}}{4}\right)r^2\\	
	&=  \frac{\sqrt{3}}{16}  \left(5-\sqrt{21}\right)r^2,	
\end{align*}
therefore

\begin{equation}\label{equ-3-g}
	S_{\triangle OED}=\frac{\sqrt{3}}{16}  \left(5-\sqrt{21}\right)r^2.
\end{equation}

It remains to compute the area of the segment $\Gamma_0$ which is a part of a circle with radius $r$ separated by the   arc $\arc{ED}$. Let $\alpha$ be the central angle of the segment $\Gamma_0$. It follows from \cref{equ-3-ea} and \cref{equ-3-f} that
\[ \sin\left(\frac{\alpha}{2}\right)=\frac{\overline{EH}}{\overline{EG}}=\frac{\left(\frac{\sqrt{14}-\sqrt{6}}{4}\right)r}{2r}=\frac{\sqrt{14}-\sqrt{6}}{8}, \]    
hence
\begin{equation}\label{equ-3-i}
	\alpha= 2\arcsin\left(\frac{\sqrt{14}-\sqrt{6}}{8}\right).
\end{equation}   
The area of the segment $\Gamma_0$ is the difference between the area of the sector $\Gamma_1$  in quadrant $GPEDL$ separated by the arc $\arc{ED}$ and that of the triangle $\triangle GED$. First note that the area of the sector $\Gamma_1$ is equal to the area of the quadrant times the ratio of their central angles, i.e., $ \frac{2\alpha}{\pi}$.   So, it follows from   \cref{equ-3-eb}, \cref{equ-3-f} and \cref{equ-3-i} that 
\begin{align*}
	S_{\Gamma_0}&= \left(\dfrac{2\alpha}{\pi}\right)\times S_{GPEDL} - S_{\triangle GED}\\
	&= \left(\dfrac{2\alpha}{\pi}\right)\times \dfrac{\pi r^2}{4} - \frac{\overline{ED}\times \overline{GH}}{2}\\	
	&=  \dfrac{\alpha r^2}{2}   -     \left(\frac{\sqrt{3}}{2} x + \dfrac{\sqrt{2}}{2}r\right)\times \left(\dfrac{x}{2}\right) \\	
	&=  \dfrac{\alpha r^2}{2}  -  \frac{\sqrt{3}}{4} x^2 - \dfrac{\sqrt{2}}{4}   xr \\
	&=  \dfrac{\alpha r^2}{2}  -  \frac{\sqrt{3}}{4} \left(\frac{\sqrt{14}-\sqrt{6}}{4}\right)^2 r^2 - \dfrac{\sqrt{2}}{4}   \left(\frac{\sqrt{14}-\sqrt{6}}{4}\right)r^2 \\
	&=  \dfrac{\alpha r^2}{2}  -  \frac{\sqrt{3}}{4} \left(\frac{5-\sqrt{21}}{4}\right)  r^2 -   \left(\frac{\sqrt{7}-\sqrt{3}}{8}\right)r^2 \\
	&=  \dfrac{\alpha r^2}{2}  -    \left(\frac{5\sqrt{3}-3\sqrt{7}}{16}\right)  r^2 -   \left(\frac{\sqrt{7}-\sqrt{3}}{8}\right)r^2 \\
	&=  \dfrac{\alpha r^2}{2}  +    \left(\frac{\sqrt{7}-3\sqrt{3}}{16}\right)  r^2.		
\end{align*}
Therefore, we get
\begin{equation}\label{equ-3-j}
	S_{\Gamma_0}= \dfrac{\alpha r^2}{2}  +   \left(\frac{\sqrt{7}-3\sqrt{3}}{16}\right)  r^2.	
\end{equation} 
If we denote the circular part consisting of $\Gamma_3$ and  $\Gamma_0$ by $\Lambda_1$, then  it follows  from \cref{equ-3-g}  and \cref{equ-3-j} that
\[ 	S_{\Lambda_1}= \dfrac{\alpha r^2}{2}  +   \left(\frac{\sqrt{7}-3\sqrt{3}}{16}\right)  r^2+ \left(\dfrac{5\sqrt{3}-3\sqrt{7}}{16}\right)r^2, \]
or equivalently
\begin{equation}\label{equ-3-k}
	S_{\Lambda_1}= \left[\arcsin\left(\frac{\sqrt{14}-\sqrt{6}}{8}\right) +   \frac{\sqrt{3}-\sqrt{7}}{8}\right] r^2.	
\end{equation}   
Now, we can use \cref{equ-3-k} to compute the area of   regular  convex 6-gon $\Lambda_6$. Since $S_{\Lambda_6} =6\times S_{\Lambda_1}$, we obtain
\begin{equation}\label{equ-3-l}
	S_{\Lambda_6}= \left[6\arcsin\left(\frac{\sqrt{14}-\sqrt{6}}{8}\right)  +   \frac{3}{4}\left(\sqrt{3}-\sqrt{7}\right)\right]   r^2.	
\end{equation}   
This formula is obtained from the general formula \cref{equ-f} if we set $n=6$.  Similarly, the area of the inscribed regular hexagon $\Gamma_6$ is computed from \cref{equ-3-g} as follows:
\begin{equation}\label{equ-3-m}
	S_{\Gamma_6}= \frac{3\sqrt{3}}{8}  \left(5-\sqrt{21}\right)r^2. 	
\end{equation}  

\begin{remark}
	It should be noted that the general formula \cref{equ-f} can be obtained using a similar argument to the one   just discussed for the case $n=6$. Another method might be the application of definite integrals. 
\end{remark}

Now, it is   time to analyze the number  given    in line 6 by the scribe as the constant of the  regular convex polyarc $\Lambda_6$. As usual, the scribe has assumed the radius $r$ to be 1 and approximated the area of $\Lambda_6$ to that of   the inscribed regular hexagon $\Gamma_6$. It should be noted that to compute the area of   $\Lambda_6$ one must  compute  the value of angle $\alpha$ which was impossible for the Susa scribe.  So the Susa scribe of this tablet has had to ignore the areas of the six circular segments in his calculations.  As we can see from \cref{equ-3-m}, there are two irrational square roots in the   formula  for the area of the inscribed regular  hexagon $\Gamma_6$. This must have made  the scribe to apply the  approximate values  in his calculations.   It is better to consider different possible approximations to  these two square roots.   

 Although    the most common approximation $ \sqrt{3}\approx \frac{7}{4}$ must have been applied here,   we  need to  check different approximations for the square root of 21. To do so, we   apply the usual Babylonian method for approximating square roots     using the formula $x_{n+1}=\frac{1}{2}(x_n+\frac{21}{x_n})$ when $x_0$ is known.\footnote{For more details, see \cite{FR98,Fri07-2,Hyp02,Mur99-1,Neu55,Neu69}.} Since $4<\sqrt{21}<5$, we may use $x_0=4,5,4\frac{1}{2}$. The first three terms  for each choice of $x_0$   are shown     in \cref{tab-SMT3-a}. Note that the decimal fractions are written to five sexagesimal places. Also, the exact value of $\sqrt{21}$ is $ 4;34,57,16,21,0,23,51,\ldots $.

\begin{table}[H] 
	\centering
	\begin{tabular}{|c|c|c|c|}
		\hline 
		& $x_1 $ & $x_2$ & $x_3 $    \\
		\hline 
		$x_0= 4$& $\frac{37}{8}\approx4;37,30 $ & $ \frac{2713}{592}\approx4;34,57,58,22,42$ & $ \frac{14720113}{3212192}\approx4;34,57,16,21,3 $    \\
		\hline
		$x_0= 5$& $\frac{23}{5}=4;36$ & $ \frac{527}{115}\approx 4;34,57,23,28,41$ & $ \frac{277727}{60605}\approx 4;34,57,16,21,0 $    \\
		\hline
		$x_0= 4;30$& $\frac{55}{12}\approx 4;35 $ & $ \frac{6049}{1320}\approx 4;34,57,16,21,49$ & $ \frac{73180801}{15969360}\approx 4;34,57,16,21,0 $       \\ 
		\hline
	\end{tabular}
	\caption{Possible approximations of $ \sqrt{21}  $}
	\label{tab-SMT3-a}
\end{table}

By using the approximations of $\sqrt{21}$ given in  \cref{tab-SMT3-a} and the usual approximation $\sqrt{3}\approx\frac{7}{4}$, we can  compute   the  possible  approximations  of $\frac{3\sqrt{3}}{8}  \left(5-\sqrt{21}\right)$ which are  shown in \cref{tab-SMT3-b}.

\begin{table}[H] 
	\centering
	\begin{tabular}{|c|c|c|c|}
		\hline 
		& $x_1 $ & $x_2$ & $x_3 $    \\
		\hline 
		$x_0=4$&	$0;14,45,56,13$& $0;16,25,18,51,4 $ & $ 0;16,26,9,53,40$     \\
		\hline
		$x_0=5$&	$0;15,45$     & $0;16,26,5,13,2 $ & $0;16,26,9,53,42 $     \\
		\hline
		$x_0=4;30$&	$0;16,24,22,30$&  $0;16,26,9,53,10$  & $0;16,26,9,53,42 $     \\ 
		\hline
	\end{tabular}
	\caption{Possible approximations of $\frac{3\sqrt{3}}{8}  \left(5-\sqrt{21}\right)$}
	\label{tab-SMT3-b}
\end{table}

It follows from the data of  \cref{tab-SMT3-b} that a reasonable approximate value of $S_{\Gamma_6} $ up to five sexagesimal places is $0;16,26,9,53$. One  might expect  the scribe to have obtained this value, however,  the scribe's value  is $0;16,26,46,40$ which is  a little bigger than the obtained approximation. Two reasonable conjectures for this error are either that (1) the scribe has used different approximations of $ \sqrt{3}$ and $\sqrt{21}$ or  (2) he  has made a mistake in his calculations.

Let us consider different approximations for two square roots $ \sqrt{3}$ and $\sqrt{21}$ and compute the corresponding approximations of  $\frac{3\sqrt{3}}{8}  \left(5-\sqrt{21}\right)$. We have considered 11 initial values  for $\sqrt{21}$ and 13 initial values for $\sqrt{3}$. The following table shows the corresponding approximate value for each square root and $\frac{3\sqrt{3}}{8}  \left(5-\sqrt{21}\right)$. As we can see from  \cref{tab-SMT3-c}, none of the obtained values equal to the scribe's value $0.274104938$, but the closest one is given in the shaded cell, i.e., $ 0.2740178571$.

\begin{table}[H] 
	\centering
	\includegraphics[scale=1]{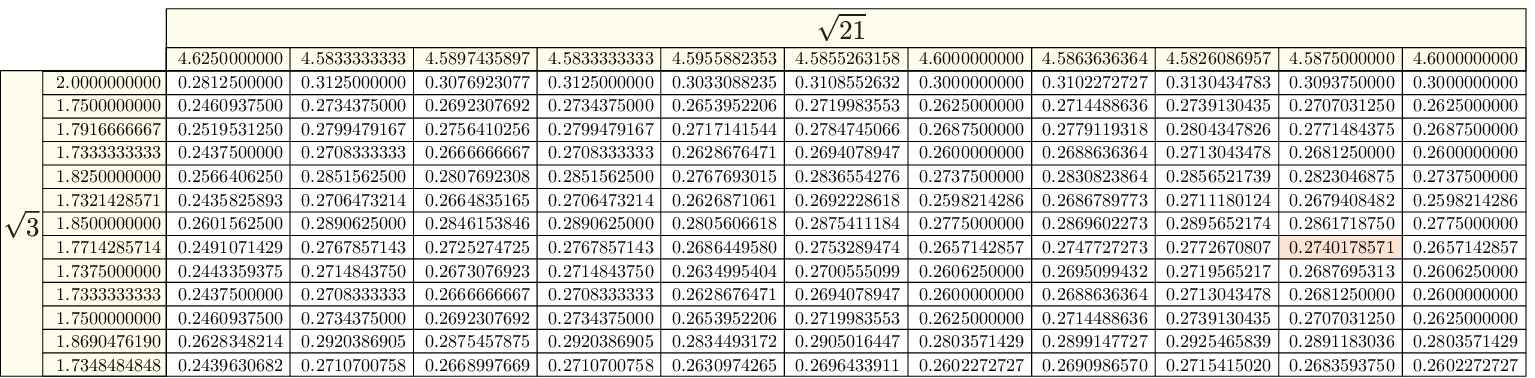}
	\caption{Different approximations of $\frac{3\sqrt{3}}{8} \left(5-\sqrt{21}\right)$} 
	\label{tab-SMT3-c}
\end{table}

From the data of   \cref{tab-SMT3-c} and the discussion above we can come to the conclusion that the scribe must have made a mistake in the process of approximating $\sqrt{21}$, because he   had to perform rather complicated calculations in the sexagesimal system.

Finally, we may compute the exact area of the  regular  convex 6-gon. If we set $r=1$ in \cref{equ-3-l}, then  $S_{\Lambda_6}=  0;17,17,12,44,20,29,36,54,\cdots $. Hence, the error percentage of scribe's approximation is
\[ e=\left|\frac{0;17,17,12,44,20,29,36,54,\cdots -0;16,26,46,40}{0;17,17,12,44,20,29,36,54,\cdots}\right|\times 100\%\approx 1.4\% \] 
which, despite the scribe's mistake,  is a relatively good approximation.

\subsection{An \textit{apusamikkum}: a hole of lyre}
The transliteration and translation of   lines  22-24 are as follows:
\begin{note1} 
	\underline{Obverse: Lines 22-24} \\
	Transliteration:\\
(L22) 26,40 igi-gub \textit{\v{s}\`{a} a-pu-s\`{a}-am-mi-ki}\\
(L23) 1,20 bar-d\'{a} \textit{\v{s}\`{a} a-pu-s\`{a}-mi-ki}\\
(L24) 33,20 \textit{pi-ir-ku \v{s}\`{a} a-pu-s\`{a}-mi-ki}
\end{note1}
\noindent
Translation:\\
(L22) 0;26,40 is the constant of (a figure called) \textit{apusamikkum}.\\
(L23) 1;20 is the diagonal  of (a figure called) \textit{apusamikkum}.\\
(L24) 0;33,20 is the transversal of (a figure called) \textit{apusamikkum}.\\

\begin{figure}[H]
	\centering
	\includegraphics[scale=1]{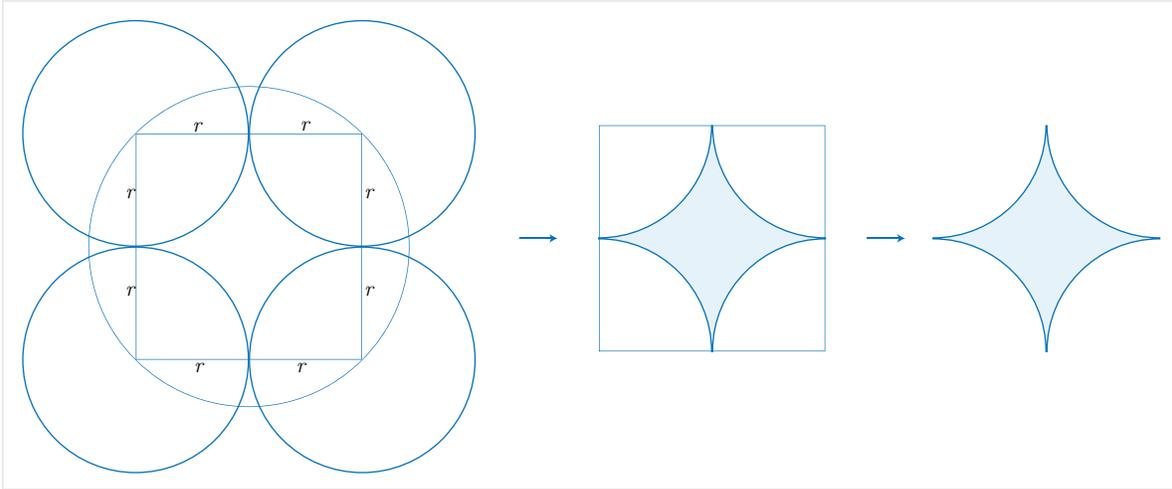}
	\caption{Construction of a  regular concave 4-arc: \textit{apusamikkum}}
	\label{Figure10}
\end{figure}

  The figure  described in these lines is called an \textit{apusamikkum}, which is a  regular  concave 4-arc. It is constructed by repeating (rotating) a circle  of radius $r$ four times around a bigger circle   such that they are pairwise tangent. It also can be made by removing four equal quarters of radius $r$ from a square of side $2r$ (see \cref{Figure10}).  It should be noted that the Akkadian word  \textit{apusamikkum}  is derived from the Sumerian word  {\fontfamily{qpl}\selectfont\'{a}b-z\`{a}-m\'{i}} ``the hole of a lyre'', where {\fontfamily{qpl}\selectfont\'{a}b $ \approx $ ab} means ``hole'' (see \cite{Mur13}).

\begin{figure}[H]
	\centering
	\includegraphics[scale=1]{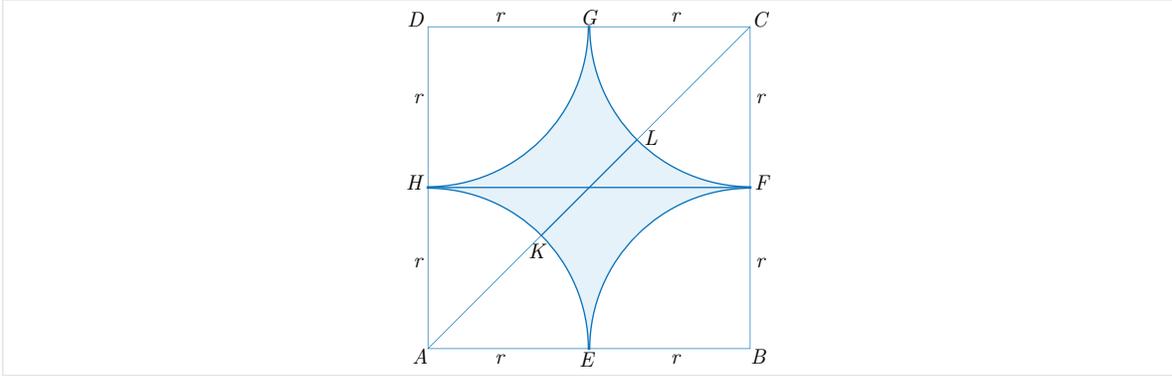}
	\caption{Area of an \textit{apusamikkum}}
	\label{Figure11}
\end{figure}

 In lines 22-24, three numbers are listed as the  area, the diagonal and the transversal of the   \textit{apusamikkum} $\Phi$. We computed the area for the general case in Section 2. If we set $n=4$ in \cref{equ-e} and simplify, we get
\begin{equation}\label{equ-3-gee}
	S_{\Phi}=\frac{(16-4\pi)a^2}{\pi^2}
\end{equation}
where $a$ is the arc length of each quarter. The other two numbers are called the diagonal and the transversal of the polyarc $\Phi$. The  transversal seems to be the intersection of the diagonal   of the circumscribed square with the  \textit{apusamikkum} $\Phi$, which is the line segment $KL$ in \cref{Figure11}. The diagonal is the horizontal or vertical diagonals $HF$  or $EG$ whose length is $2r$. So, 
\begin{equation}\label{equ-3-ge}
	\overline{HF}=2r.	
\end{equation}

To compute the transversal $KL$, first note that 
$\overline{KL}=\overline{AC}-(\overline{AK}+\overline{CL}), $
and since $\overline{AK}=\overline{CL}=r$, we get
\begin{equation}\label{equ-3-gea}
	\overline{KL}=\overline{AC}-2r.	
\end{equation}
It remains to find $\overline{AC}$.  By the Pythagorean theorem in the right triangle $ \triangle ABC$, since $\overline{AB}=\overline{CB}=2r $,   we have
\begin{equation}\label{equ-3-gf}
	\overline{AC}=2\sqrt{2}r.	
\end{equation}
Now, since  $a$ is the quarter of the circumference of a circle with radius $r$, we have$a=\frac{\pi r}{2}$ or $r=\frac{2a}{\pi}$. So, \cref{equ-3-gea} and \cref{equ-3-gf} imply that
\begin{align*}
	\overline{KL} &=\overline{AC}-2r\\
	&=  2\sqrt{2}r - 2r \\
	& = 2(\sqrt{2} -1)r \\
	&=  2(\sqrt{2} -1)\times \frac{2a}{\pi},  
\end{align*}
or
\begin{equation}\label{equ-3-geb}
	\overline{KL}= \dfrac{4(\sqrt{2}-1)a}{\pi}.	
\end{equation}

Now, we can check the scribe's numbers by assuming $a=1$ and using the Babylonian approximations $\pi\approx 3$ and $\sqrt{2}\approx \frac{17}{12}=1;25$. It follows from \cref{equ-3-ge} that
\[ \overline{HF}=\dfrac{4a}{\pi}\approx\frac{4\times 1}{3}=4\times (0;20)=1;20 \]
which is the number in line 23. By \cref{equ-3-geb}, we have
\[ 	\overline{KL}= \dfrac{4(\sqrt{2}-1)a}{\pi}\approx   \dfrac{4\times (1;25-1)\times 1}{3}=\frac{4\times (0;25)}{3}=(1;20)\times (0;25)=0;33,20 \] 
which is the number in line 24. The approximate area is obtained from \cref{equ-3-gee}:
\[ S_{\Phi}=\frac{(16-4\pi)a^2}{\pi^2} \approx \frac{(16-12)\times 1}{9}=\frac{4}{9}= 0;26,40 \]
which is listed in line 22. Again, the numbers of the Susa scribe here are correct approximate values of the  area, the diagonal and the transversal of an   \textit{apusamikkum} whose arc length is unit.

\subsection{An \textit{apusamikkum} with three vertices}
The transliteration and translation of   line   25 are as follows:
\begin{note1} 
	\underline{Obverse: Line 25} \\
	Transliteration:\\
(L25) 15 igi-gub \textit{\v{s}\`{a} a-pu-s\`{a}-mi-ik-ki \v{s}\`{a}} 3
\end{note1}
\noindent
Translation:\\
(L25) 0;15 is the constant of (a figure called) \textit{apusamikkum}  of three vertices.\\

\begin{figure}[H]
 	\centering
 	\includegraphics[scale=1]{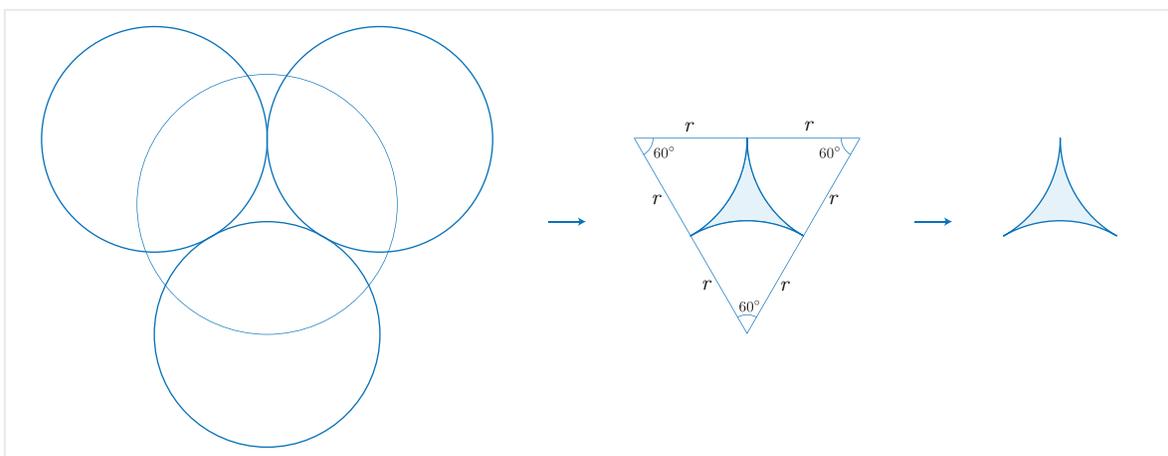}
 	\caption{Construction of an \textit{apusamikkum} with three vertices}
 	\label{Figure12}
 \end{figure}

  The figure    in this line   is a regular concave 3-arc. Beside the standard way, it can be made by removing three equal one-sixths of a circle with radius $r$ from a equilateral triangle with side $2r$ (see \cref{Figure12}).   
  
  In line 25, the area of this figure $\Omega$ is listed assuming that the radius of the circle is $r=1$. The exact value of the area is obtained from \cref{equ-e} by setting $n=3$ and $a=\frac{\pi r}{3}$:
  \begin{align*}
  	S_{\Omega}&= \left[\frac{27}{\pi^2}\times \cot\left(\frac{\pi}{3}\right)-\frac{9}{2\pi}\right]a^2\\
  	&= \left(\frac{27}{\pi^2\sqrt{3}}-\frac{9}{2\pi}\right)\times\frac{\pi^2r^2}{9} \\
  	&=  \sqrt{3}r^2 -  \frac{\pi r^2}{2},
  \end{align*}
 or 
   \begin{equation}\label{equ-3-gef}
  	S_{\Omega}=\left(\sqrt{3}  -  \frac{\pi}{2}\right)r^2.	
  \end{equation}

If we set $r=1$ in \cref{equ-3-gef} and use the Babylonian approximations  $ \sqrt{3}\approx\frac{7}{4} $  and $\pi\approx3 $,  we get
\begin{align*}
	S_{\Omega}=\left(\sqrt{3}  -  \frac{\pi}{2}\right)r^2 \approx  \left(\frac{7}{4} -  \frac{3}{2}\right)\times 1= 1;45-1;30=0;15
\end{align*}
which is the number listed in line 25   verifying the Scribe's number.

\section{Significance of circular figures} 
 All the calculations and explanations we have given in the previous sections demonstrate the mathematical importance of these constants. Every mathematician or   student of mathematics  understands the high level of mathematical skill needed to deal with these geometric figures. The Susa scribes were  the skilled ``geometers'' of their age and we believe this remarkable list of geometric constants identifies the special importance of geometry in   Elamite mathematic. 
 
 Surprisingly, some of these circular figures were used to decorate Elamite artifacts. In the book ``The Art of Elam'' published in 2020, there are  pictures of many artifacts excavated from ancient sites   (mostly Susa) in the southwest of Iran. Some of these artifacts use the circular figures that we discussed here.  The following table gives the photo reference for some of these usages. 
\begin{table}[H] 
	\centering
	\begin{tabular}{|c|c|}
		\hline 
	Circular Figure	& Photo Reference        \\
		\hline 
	barley-field 	& \cite{Alv20}:  Plate 12-c,  Plate 12-g \\
	\hline
	\textit{apusamikkum}  	&\cite{Alv20}:   Plate 12-b, Plate 49-f \\
	       with four vertices                          	&  Louvre Collection: SB23397, A58786\\
	\hline
		\textit{apusamikkum}	&\cite{Alv20}: Plate 11-e, Plate 11-j,  Plate 11-m \\
				   with three vertices                                                   	& Louvre Collection: SB19519, SB23369, SB23386, SB23462\\
	\hline
	\end{tabular}
	\caption{Examples of circular figures in Elamite artifacts}
	\label{tab-c}
\end{table}  
 
It is interesting that the decoration on some artifacts use similar constructions to those of the  circular figures we have discussed in this article. For example, inside two terracotta vases\footnote{The reference numbers of these item are SB23369 and SB23386. The reader can see their photos on the Louvre collection website: \url{https://collections.louvre.fr/ark:/53355/cl010174158} and \url{https://collections.louvre.fr/ark:/53355/cl010174136}.} dated to ca. 4000--3500 BC,  three almost pairwise tangent circles   enclose a small  regular concave 3-arc. 

\begin{figure}[H]
 	\centering
 	\includegraphics[scale=1]{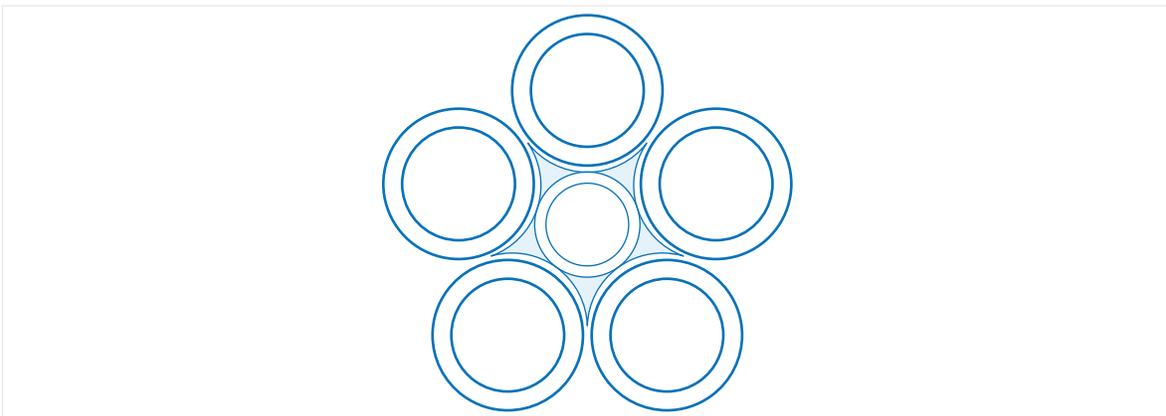}
 	\caption{Reconstruction of the pattern inside a terracotta vase: SB23397}
 	\label{SB23397}
 \end{figure}

A similar theme  is  found inside another terracotta vase\footnote{The reference number of this item is SB23397 and the reader can see its photo on the Louvre collection website: \url{https://collections.louvre.fr/ark:/53355/cl010173875}.} from the same era. In this pattern, a chain of five equal repeating circles encloses another circle and the space  between five surrounding  circles and the middle one  is occupied by    five        regular  concave 3-arcs whose union forms the corners of a     regular concave 5-arc.  A reconstruction of this pattern is shown in \cref{SB23397}. 
 
Besides these examples, there are many other Elamite artifacts whose decorations contain  different geometric patterns. The examples given  in this section may be followed up by the interested reader  either through visiting    the Louvre collection website or consulting the book ``The Art of Elam'' to see additional such patterns (see \cite{Alv20}). In an ongoing book project on Elamite mathematics, the authors have attributed two chapters to categorize the geometric patterns used in   Elamite art. This part of the book demonstrates the wide variety of geometric shapes and  patterns that Elamite artists used in their artworks implying a  close connection between geometry and Elamite art.

\section{Conclusion} 
The detailed analysis of some of the geometric constants in  this tablet  clearly suggests  that the Susa scribes  enjoyed   considerable knowledge about two-dimensional figures and   plane geometry. In particular, computing the area   of       complex geometric  figures   reveals  that the scribes of Susa  knew  those  figures very well and the methods needed to approximate their areas.  As we saw, some of these calculations are  complicated and performing them requires a high level of mathematical skill and experience.

As Bruins   rightly has  pointed out   \cite{Bru50}, the list of geometric  constants  in   this text provides a documentary evidence that Susa scribes    performed some of their calculations from a purely  mathematical point of view, or in his words ``science for science''.      This observation     suggests  that  in  ancient Susa      mathematical skills were most likely taught   in  an  educational setting   which might reasonably be named  ``The Susa School of Mathematics''\footnote{As is known from Mecquenem's report  (see the introduction of the \textbf{TMS} \cite{BR61}), all the  26 Susa mathematical  tablets   were collected in 1933 under a large paving $30m$ wide and $15m$ long, at the  Ville Royale, where, in 1934--35, father Van der Meer also collected some school texts which he published in  the M\a'{e}moires de la Mission Arch\a'{e}ologique de Iran: Volume XVII. All these artifacts taken together suggest the existence of a special educational establishment  in  ancient Susa.}.

Elamite artists and craftsmen  were  obviously interested in using various complicated and beautiful geometric patterns to decorate their works and this application of geometric patterns in Elamite art is strong evidence for the relationship between geometry and art in   ancient Elam.

\end{document}